\documentclass[12pt]{amsart}

\usepackage{a4}
\usepackage{amssymb}
\usepackage{amsmath}
\usepackage{amsthm}
\usepackage{amstext}
\usepackage{amscd}
\usepackage{latexsym}
\usepackage{graphics}

\swapnumbers
\theoremstyle{plain}
\newtheorem{thm}{Theorem}[section]

\newtheorem{lemma}[thm]{Lemma}
\newtheorem{lem}[thm]{Lemma}

\newtheorem{prop}[thm]{Proposition}

\theoremstyle{definition}
\newtheorem{rem}[thm]{Remark}
\newtheorem{point}[thm]{}

\newtheorem{defn}[thm]{Definition}
\newtheorem{ex}[thm]{Example}
\newtheorem{claim}[thm]{Claim}

\newcommand{\tensor}{\otimes}

\newcommand{\Union}{\bigcup}
\newcommand{\intersection}{\cap}

\newcommand{\ilim}{\mathop{\varprojlim}\limits} 
\newcommand{\hocolim}{\mathop{\underrightarrow{\rm hlim }}}  
\newcommand{\holim}{\mathop{\underleftarrow{\rm hlim }}}  

\newcommand{\Hom}{{\rm Hom}}

\newcommand{\Spec}{{\rm Spec \,}}

\renewcommand{\tilde}{\widetilde}
\newcommand{\sA}{{\mathcal A}}
\newcommand{\sB}{{\mathcal B}}
\newcommand{\sC}{{\mathcal C}}

\newcommand{\sE}{{\mathcal E}}
\newcommand{\sF}{{\mathcal F}}
\newcommand{\sG}{{\mathcal G}}

\newcommand{\sM}{{\mathcal M}}

\newcommand{\sO}{{\mathcal O}}

\newcommand{\sT}{{\mathcal T}}

\newcommand{\N}{{\mathbb N}}

\newcommand{\Z}{{\mathbb Z}}

\input{xy}
\xyoption{all}

\begin{document}
\title{Products, Homotopy limits and Applications}
\author{Amit Hogadi  \ \ and \ \ Chenyang Xu}
\begin{abstract}In this note, we discuss the derived functors of infinite
products and homotopy limits. $QC(X)$, the category of quasi-coherent sheaves on a Deligne-Mumford stack $X$, usually has the property that the derived functors of product vanish after a finite stage. We use this fact to study the convergence of certain homotopy limits and apply it compare the derived category of $QC(X)$ with certain other closely related triangulated categories. 
\end{abstract}
\maketitle

\tableofcontents

\section{Introduction}

\noindent Let $\sA$ be an abelian category with products and enough
injectives and $D(\sA)$ be its unbounded derived category. One of
the tools used to study the unbounded derived category is the concept of
a homotopy limit as defined in \cite{bneeman}. If the
abelian category $\sA$ satisfies AB4*, i.e. if products are exact in
$\sA$, then it is trivial to observe that for any object $N\in
D(\sA)$, a map
$$ N\to \holim_{i\geq 0} N_{\geq -i}$$ is a quasi-isomorphism. This is useful because it allows us to reduce certain questions related to unbounded complexes to bounded below complexes. \\

\noindent For an algebraic stack $X$, we will denote by $QC(X)$ the
category of quasi-coherent sheaves on $X$. Unfortunately, $QC(X)$
rarely satisfy AB4*. However, in several interesting examples, as
shown by the following theorem, these abelian categories satisfy a
weaker axiom, AB4*-$n$, for some $n\geq0$. This just means that the
$i$-th derived functors of the product vanish if $i>n$ (see Section $2$).

\begin{thm}\label{thm:dm-ab4*n}
Let $S$ be a separated noetherian scheme and $X/S$ be a separated
Deligne-Mumford stack of finite type. Assume
\begin{enumerate}
 \item The coarse moduli space of $X$ is a scheme.
 \item $X$ is covered by Zariski open substacks each of which admits a finite \'etale cover from a scheme.
\end{enumerate}
Then $QC(X)$ satisfies AB4*-$n$ for some positive integer $n$.
\end{thm}

\begin{rem}
 In particular, if $S=\Spec(k)$ is a field, the conditions of the above theorem are satisfied for a quotient stack $X/k$ which is a separated Deligne-Mumford stack with quasi-projective coarse moduli space (see \cite{kresch}, Proposition 5.2). Recall that a stack $X$ is called a {\it  quotient stack} if it can be expressed as a quotient of a scheme by a linear algebraic group. In the case when the stack is noetherian and normal, this is equivalent to saying that $X$ has resolution property, i.e., every coherent sheaf can be expressed as quotient of a locally free sheaf (see \cite{totaro}).
\end{rem}

\noindent In the following theorem, we observe that the weaker
condition AB4*-$n$, is still enough to make the concept of homotopy
limits useful for reducing questions about unbounded complexes to
bounded below complexes.

\begin{thm}\label{thm:ab4*-n}
Let $\sA$ be an abelian category with enough injectives. Assume
$\sA$ satisfies AB4*-$n$ for some positive integer $n$. Then for
every $N\in D(\sA)$, we have an isomorphism $N\to \holim N_{\geq -i}$ in
$D(\sA)$. In particular, the category of $\sA$-complexes has enough $K$-injective objects. (For
the definition of $K$-injective objects, see \cite{spaltenstein})
\end{thm}

\noindent For an algebraic stack $X$ we denote by $D(X)$ the
unbounded derived category of quasi-coherent sheaves on $X$.
$D(\sO_X)$ (resp. $D^{cart}(\sO_X)$) will stand for the unbounded
derived category of (resp. cartesian) $\sO_X$-modules on the
lisse-\'etale site of $X$. Let $D_{qc}(\sO_X)$ (resp.
$D^{cart}_{qc}(\sO_X)$) be the full subcategory of $D(\sO_X)$
$D^{cart}(\sO_X)$ consisting of complexes with cartesian
quasi-coherent cohomology. We apply Theorem (\ref{thm:ab4*-n}) to
prove the following result, which is the main result of this paper. 

\begin{thm}\label{thm:comparison}
$X$ be any algebraic stack such that the diagonal morphism
$\Delta_{X/\Z}:X\to X\times_{\Z}X$ is affine. Assume that $QC(X)$
satisfies AB4*-$n$ for some positive integer $n$. Then we have a
natural equivalence
$$ D(X) \cong D_{qc}(\sO_X)\cong D^{cart}_{qc}(\sO_X).$$
\end{thm}

\noindent The above theorem is already known in the case separated quasi-compact schemes (cf. \cite{bneeman}). To prove the above theorem we introduce the notion of a {\it strong Serre subcategory} (see Section \ref{section:sssub}). The diagonal morphism of a stack being affine will actually imply
that $QC(X)$ is a strong Serre subcategory of categories of
cartesian $\sO_X$-modules and general $\sO_X$-modules on the
lisse-\'etale site of $X$. Theorem (\ref{thm:comparison}) is then
the result of the following general comparison theorem.
\begin{thm}\label{thm:sssub}
Let $F:\sC \to \sM$ be an embedding of a strong Serre subcategory of $\sM$. Assume $F$ has a right adjoint. Let $D_{\sC}(\sM)$ denote the full subcategory of $\sM$ consisting of complexes whose cohomology lies in $\sC$. Assume both $\sC$ and $\sM$ satisfy AB4*-$n$ for some nonnegative integer $n$. Then the natural functor $$ RF: D(\sC) \to D_{\sC}(\sM)$$
is an equivalence.
\end{thm}

\noindent {\bf Acknowledgement}: We thank Bhargav Bhatt and Jan-Erik Roos for their useful comments and emails. Parts of this work were done during our stay in University of Utah and Mathematical Sciences Research Institute (MSRI). We thank these institutes for their hospitality. One of our motivations to start this project was to prove that the derived category of a Deligne-Mumford stack under suitable assumptions is compactly generated. We have been told by Amalendu Krishna that he has proved this. Our paper may have some overlap with his preprint.

\section{Preliminaries}
\noindent In this section we recall some well known facts about the
derived functors of product and homotopy limits. For an abelian
category $\sA$ we use the notations $K(\sA)$ and $D(\sA)$ to denote
the category of complexes up to homotopy and the derived category
respectively. A superscript `$+$' (e.g. $D^+(\sA)$) will denote the
respective bounded below version. \\

\noindent {\bf Preliminaries on Derived Functors of Product} 
\begin{point}[Derived Functors of Product]\label{p:fpv}
Let $\sA$ be an abelian category with enough injectives. Assume small products
exist in $\sA$. For any indexing set $I$,
let $\sA^I$ be the $I$-fold product of $\sA$. $\sA^I$ itself is an abelian category
with enough injectives and since products exist in $\sA$, we have a natural functor
$$\Pi:\sA^I \to \sA \hspace{3mm} \text{given by} \ \ (A_i)_{i\in I} \to \mathop{\Pi}_{i}A_i$$
This functor is additive and since it admits a right adjoint, it is left exact.
We will denote by $\Pi^{(n)}A_i$ the $n$-th right derived functor of the above functor
evaluated at $(A_i)_{i\in I}$. Concretely, in order to calculate $\Pi^{(n)}A_i$ one chooses
injective resolution $0 \to A_i \to I^{\bullet}_i$ of $A_i$ for each $i$ and then
$\Pi^{(n)}A_i= H^n(\mathop{\Pi}_i I_i^{\bullet})$.  \\
\end{point}
\begin{point}\label{p:2ss}
Let $F:\sA \to \sB$ be an additive functor between abelian categories.  Assume $\sA$ has small products and enough injectives. Let $\{A_{\alpha}\}$ be a set of  objects in $\sA$ and for each $A_{\alpha}$ and $$ 0 \to A_{\alpha} \to C^{\bullet}_{\alpha} $$ be a resolution. Then we have the following spectral sequence
$$ E^{pq}_1=\Pi^{(q)}(C^p)\Rightarrow \Pi^{(p+q)}A_{\alpha}$$

\end{point}

\begin{defn}\cite{janerik}
Let $\sA$ be an abelian category with enough injectives and in which products
exist. For a nonnegative integer $n$,
$\sA$ is said to satisfy the {\it axiom  AB4*-$n$} if
for any set $I$, and any collection of objects $\{A_{\alpha}\}_{\alpha \in I}$
of $\sA$
$${\mathop{\Pi}_{\alpha}}^{(q)}A_{\alpha}=0 \hspace{5mm} \forall \ q\geq N+1$$
\end{defn}

\begin{rem}Let $A$ be a noetherian local ring and let $X$ be the complement of the closed point in $\Spec(A)$.
Roos showed that $QC(X)$ is AB4*-$n$ (cf. \cite{janerik}, Theorem 1.15), where $n=\max\{\dim(A),1\}$, and this bound is sharp.
\end{rem}

\noindent {\bf Preliminaries on Homotopy Limits}
\begin{defn}\cite{bneeman}\label{def:hlim}
Let $\sT$ be a triangulated category. Let
$$\cdots \stackrel{\alpha}{\to} M_{i+1}\stackrel{\alpha}{\to} M_i \stackrel{\alpha}{\to} \cdots \stackrel{\alpha}{\to} M_2 \stackrel{\alpha}{\to} M_1 $$
be an inverse system of objects $M_i \in \sT$ indexed by the natural numbers
$\N$. Assume $\Pi \ M_i$ is representable in $\sT$. Then the {\it homotopy limit} of $\{M_i\}$,
denoted by $\holim M_i$ is defined by the following exact traingle in $\sT$.
$$  \holim M_i \to  \Pi \ M_i \stackrel{1-\alpha}{\longrightarrow} \Pi \ M_i \to \holim M_i [1] $$
Dually one can also define the notion of a homotopy colimit of
a sequence $M_1\to M_2 \to \cdots $ by the following distinguished triangle.
$$ \oplus M_i \stackrel{1-\alpha}{\longrightarrow} \oplus M_i \to \hocolim M_i \to \oplus M_i[1]$$
\end{defn}

\begin{rem}\label{p:bneeman}
We recall the following features of homotopy limits from \cite{bneeman}.
\begin{enumerate}
\item Since the cone of any morphism in a traingulated category is only unique
up to a  non-canonical isomorphism, homotopy limits are only unique up to a {\it non-canonical} isomorphism.
\item Let $\{M_i\}_{i\in \N}$ be an inverse system in $\sT$ such that $\holim M_i$
exists. Let $\phi_i: L \to M_i$ be a collection of morphisms from $L$ to
$M_i's$. $\phi_i$ is a compatible system of morphisms iff the composition
$\Pi \phi: L \to \Pi \ M_i \stackrel{1-\alpha}{\longrightarrow} \Pi \ M_i$ is zero.
Thus in this case there is an induced morphism (non-unique in general)
$\phi : L \to \holim M_i$. The non-uniqueness of $\phi$ means that in general
the natural map $\Hom_{\sT}(L,\holim M_i) \to \ilim \Hom_{\sT}(L,M_i)$ is {\it not}
an isomorphism. In other words homotopy limit is usually not a limit. However whenever $RHom$ makes sense,
$RHom_{\sT}(L,\holim M_i) \cong \holim RHom_{\sT}(L,M_i)$.
\end{enumerate}
\end{rem}

\begin{rem}\label{bbprod}
Let $\sA$ be an abelian category with products. Let $\{M_i\}$ be a
collection of complexes in $K(\sA)$ such that each $M_i$ is
$K$-injective. Then it is easy to see that the class of the product
complex $\Pi M_i$ in $D(\sA)$ represents the product of $M_i's$ in
$D(\sA)$. In particular, if each $M_i$ is bounded below, and $\sA$
has enough injectives, then $\Pi M_i$ exists in $D(\sA)$. Thus for
any inverse system $\{M_i\}_{i\in \N}$ of objects in $D(\sA)$,
$\holim M_i$ exists if each $M_i$ is bounded below. One can also talk about homotopy limits in $K(\sA)$ and one can show that the homotopy limit of $K$-injective objects in $K(\sA)$ is again a $K$-injective.
\end{rem}


\section{Products of Quasi-coherent sheaves on a Deligne-Mumford stack}

\noindent Throughout this section we fix a noetherian separated base scheme $S$ and let $X/S$ be a separated Deligne-Mumford stack of finite type. We denote by $q:X\to \underline{X}$ the coarse moduli space of $X$. In particular, $\underline{X}$ is separated.

\begin{lem}\label{lem:affine}
Let $i:U=q^{-1}(\underline{U})\to Y$ be an open substack such that the coarse moduli space $\underline{U}$ of $U$ is an affine open subset of $\underline{Y}$. Then the inclusion $i:U\to Y$ is affine. In particular, the functor $$i_*:QC(U)\to QC(Y)$$ is exact.
\end{lem}
\begin{proof}Left to the reader.
\end{proof}

\begin{lem}\label{lem:cover}
Let $Y/S$ be any stack and $Y=\Union_{i=1}^rV_i$ be a Zariski cover by finitely many open substacks. For $1\leq k \leq r$, let
$$U_k = \mathop{\amalg}_{1\leq i_1<i_2 ...<i_k\leq r} V_{i_1}\intersection \cdots V_{i_k}$$ and $j_k:U_k\to Y$ be the natural maps.
 Assume the following holds.
\begin{enumerate}
 \item[(i)] For every $k$, $j_{k*}:QC(U_k)\to QC(Y)$ is exact.
 \item[(ii)] For every $k$, $QC(U_k)$ is AB4*-$n_k$ for some positive integer $n_k$.
\end{enumerate}
Then $QC(Y)$ is AB4*-$n$ for any $n\ge\max_k\{k+n_k\}$.
\end{lem}
\begin{proof}
Let $\{F_{\alpha}\}_{\alpha \in I}$ be a set of quasi-coherent
sheaves on $Y$. For each $F_{\alpha}$ we have a C\v{e}ch resolution
$$ 0 \to F_{\alpha} \to \sC^1(F_{\alpha}) \to \cdots \to \sC^r(F_{\alpha})\to 0 $$
where for each $k$, $\sC^k(F_{\alpha}) = j_{k*}j_k^*\sF$.\\

\noindent Since $j_{k*}$ is exact and has an exact left adjoint, it
preserves products and also maps injective objects in $QC(U_k)$ to
injective objects in $QC(Y)$. Since $j_{k*}$ is exact, we have
$$j_{k*}(\Pi^{(i)}j_k^*\sF)=(\Pi^{(i)})\sC^k(F_{\alpha}).$$
Thus from the assumption that $QC(U_k)$ is AB4*-$n_k$, we get
$$\Pi^{(i)}\sC^k(F_{\alpha})=0 \ \forall \ i>n_k$$

\noindent Now from (\ref{p:2ss}), we have a spectral sequence
$$E_1^{p,q}= \Pi^{(q)}\sC^p(F_{\alpha}) \Rightarrow \Pi^{(p+q)}F_{\alpha}$$
This proves the lemma.
\end{proof}

\begin{rem}\label{schemecase} If $U$ is an affine scheme, then $QC(U)$ is AB4*. Thus the above lemma immediately implies that if $Y$ is a separated quasi-compact scheme, then $QC(Y)$ satisfies AB4*-$n$, where $n+1$ is the minimum number of affine opens required to cover $X$.
\end{rem}

\begin{lem}\label{lem:globalcase}
Let $Y$ be any stack which admits a finite \'etale cover $f:Z\to Y$
where $Z$ is a separated quasi-compact scheme. Then $QC(Y)$
satisfies AB4*-$n$ for some positive integer $n$.
\end{lem}
\begin{proof}
Since $f$ is finite \'etale, $f^*$ is both right as well as left adjoint of $f_*$. Moreover both $f_*,f^*$ are exact. Thus $f^*$ preserves products and takes injectives to injectives. Therefore for any set of sheaves $\{\sF_{\alpha}\}$ in $QC(Y)$
$$ \Pi^{(i)}f^*\sF_{\alpha} = f^* \left(\Pi^{(i)}\sF_{\alpha}\right)$$
Since $QC(Z)$ satisfies AB4*-$n$ for some positive integer $n$ (\ref{schemecase}), $ \Pi^{(i)}f^*\sF_{\alpha}=0 \ \forall \ i>n$. This implies $\Pi^{(i)}\sF_{\alpha}=0 \ \forall \ i>n$.
\end{proof}

\begin{proof}[Proof of Theorem \ref{thm:dm-ab4*n}] The statement now is straightforward from the assumptions, Lemma (\ref{lem:affine}) and Lemma (\ref{lem:globalcase}). 
\end{proof}

\section{Homotopy Limit of Truncations}

\noindent In this section we prove Theorem (\ref{thm:ab4*-n}).

\begin{lem}\label{lem1}
Let $\sA$ be any abelian cateogry with products. Let $N \in D^{+}(\sA)$. Then for any integers $i_0$ and $j_0$
$$\holim_{i\geq i_0}N_{\geq -i} \cong \holim_{i\geq j_0}N_{\geq -i}$$
\end{lem}
\begin{proof}
Without loss of generality we may assume $j_0=i_0+1$. For simplicity of notation, let $N_i=N_{\geq -i}$. We now have the following diagram where rows and last two columns are distinguished triangles.
$$\xymatrix{
            0\ar[r]                  & N_{i_0+1}\ar[r]^{\cong}\ar[d]        &                        N_{i_0+1}\ar[d] \\
\holim_{i\geq i_0+1}N_i \ar[r] \ar[d] & \Pi_{i\geq i_0+1} N_i \ar[r]^{1-{\rm shift}}\ar[d] & \Pi_{i\geq i_0+1} N_i\ar[d] \\
\holim_{i\geq i_0}N_i \ar[r]        & \Pi_{i\geq i_0} N_i \ar[r]^{1-{\rm shift}}       & \Pi_{i\geq i_0} N_i
}
$$
By the octahedron axiom, the triangle
$$ 0 \to \holim_{i\geq i_0+1}N_i \to \holim_{i\geq i_0}N_i $$
must be distinguished. This proves the lemma.
\end{proof}

\begin{lem}\label{lem2}
Let $\sA$ be any abelian category with products and $A$ be an object of $\sA$. Then the sequence
$$ A \to \Pi_{i\geq 0}A \stackrel{1-{\rm shift}}{\longrightarrow} \Pi_{i \geq 0}A $$
is exact.
\end{lem}
\begin{proof}
 Left to the reader.
\end{proof}

\begin{lem}\label{bbelow}
Let $\sA$ be any abelian category and $N\in D^+(\sA)$. Then for any integer $i_0$,
$$N \cong \holim_{i\geq i_0}N_{\geq -i}$$
\end{lem}
\begin{proof}
By Lemma (\ref{lem1}) we may assume $$H^{-i}(N)=0 \ \ \forall \ i\geq i_0$$
In this case choose a complex $C^{\bullet}$ representing $N$ such that $C^{-i}=0 \ \forall i \geq i_0$. Thus $C^{\bullet}_{\geq -i}=C^{\bullet} \ \forall i \geq i_0$. The lemma now follows because by Lemma (\ref{lem2})
the following is an exact sequence of complexes
$$ 0 \to C^{\bullet} \to \Pi_{i\geq i_0}C^{\bullet}_{\geq -i} \stackrel{1-{\rm shift}}{\longrightarrow} \Pi_i C^{\bullet}_{\geq -i}\to 0$$
\end{proof}

\begin{lemma}\label{cor:cantruncate}
Let $\sA$ be an abelian category which is AB4*-$n$. Let $\{M_{\alpha}\}$ and $\{L_{\alpha}\}$ be a collection of bounded below complexes of injective objects. Let $t$ be an integer and let $\phi_{\alpha}:M_{\alpha}\to L_{\alpha}$ be a collection of morphisms such that the induced map $$H^i(M_{\alpha}) \to H^i(L_{\alpha})$$ is an isomorphism for all $i \geq t$. Then for each $i \geq t + n + 1 $, $$H^i(\Pi M_\alpha) \to H^i(\Pi L_{\alpha})$$ is an isomorphism.
\end{lemma}
\begin{proof}
For a collection of bounded below complexes $\{N_{\alpha}\}$ of injective objects, we have the following spectral sequence
$$ E_2^{p,q}=\Pi^{(p)}_{\alpha}H^q(N_{\alpha}) \Rightarrow H^{p+q}(\Pi_{\alpha} N_{\alpha}).$$
Since $E_2^{p,q}=0$ for $p<0$ or $p>n$, the spectral sequence is
convergent. Now we apply it two both $\{M_{\alpha}\}$ and $\{L_{\alpha}\}$.
\end{proof}

\begin{proof}[Proof of (\ref{thm:ab4*-n})]
Let $N\in Kom(\sA)$ be  a  complex. Fix an integer $t$, Let $N_{\geq -t}$ be the truncation complex and $L$ its injective resolution with $\tau:N\to L$ the induced map. We denote by $M_i$ an injective resolution of $N_{\ge -i}$  and let $$\tilde{\tau}:\Pi_iM_i \to \Pi_iL_{\geq -i}$$ be a map induced by $\tau$. We have the following commutative diagram.

$$ \xymatrix{
 \holim M_i \ar[r]\ar[d]^{\tau} & \Pi_i M_i \ar[r]^{1-shift}\ar[d]^{\tilde{\tau}} & \Pi_iM_i\ar[d]^{\tilde{\tau}} \\
  L \ar[r]^{h_L}                            & \Pi_i L_{\geq -i} \ar[r]^{1-shift} & \Pi_i L_{\geq -i}
}
$$

\noindent Note that the top row defines a distinguished triangle by
definition of $\holim M_i$ and the bottom row is an exact sequence 
by Lemma (\ref{bbelow}) and hence also defines a distinguished traingle. We also have maps $h_N:N\to \holim
M_i$ and $\tau:N\to L$. Thus it gives a a following diagram for each
$k$.

$$ \xymatrix{
H^k(N)\ar[rd] \ar[rdd] & &  & \\
& H^k(\holim M_i) \ar[r]^{h_M}\ar[d]^{\tau} & H^k(\Pi_i M_i) \ar[r]^{1-shift}\ar[d]^{\tilde{\tau}} & H^k(\Pi_iM_i)\ar[d]^{\tilde{\tau}} \\
& H^k(L) \ar[r]^{h_L}                            & H^k(\Pi_i L_{\geq -i}) \ar[r]^{1-shift} & H^k(\Pi_i L_{\geq -i})
}
$$

\noindent It follows from Lemma (\ref{cor:cantruncate}) that the map $H^k(\Pi_i M_i) \to H^k(\Pi_i L_{\geq -i})$ is an isomorphism for $k\ge n-t+1 $. Therefore, $H^k(\holim M_i) \cong H^k(L)$ for $k \ge n-t+2$. Since $H^k(N)\to H^k(L)$ is an isomorphism for $k\geq -t$, we see that for all $k$ large enough as compared to $-t$, $H^k(N)\to H^k(\holim M_i)$ is an isomorphism. But $t$ was arbitrary. Therefore $N\to \holim M_i$ is a quasi-isomorphism.
\end{proof}

\section{Strong Serre subcategories}\label{section:sssub}
\noindent Let $\sM$ denote an abelian category and $\sC$ denote a full subcategory of $\sM$. Recall the following definition of a Serre category.
\begin{defn}[Serre subcategory]
$\sC$ is called a {\it Serre subcategory} if for any objects $a,b$
of $\sC$, and any exact sequence
$$ 0 \to a \to c \to b \to 0$$
$c$ is also an object of $\sC$.
\end{defn}
\noindent Throughout this section we make the following assumptions on $\sM$ and $\sC$
\begin{enumerate}
\item[$\bullet$] $\sM$ and $\sC$ have enough injectives.
\item[$\bullet$] The inclusion functor $F:\sC \to \sM$ is exact.
\end{enumerate}
In this case, saying $\sC$ a Serre subcategory is equivalent to saying that for any objects $a,b$ of $\sC$, the natural map
$$ Ext^1_\sC(a,b) \to Ext^1_{\sM}(a,b)$$
is an isomorphism. This as well as examples in algebraic geometry which will be mentioned below motivate the following definition.
\begin{defn}
$\sC$ is called a {\it strong Serre subcategory} if for any two
objects $a,b$ of $\sC$, the natural maps
$$ Ext^i_{\sC}(a,b) \to Ext^i_{\sM}(a,b)$$
are isomorphisms for all $i\geq 0$.
\end{defn}
\noindent The motivation for the definition comes from the following easy proposition.
\begin{prop}\label{gsiff}
Let $\sC$ be a full subcategory of $\sM$. Then $\sC$ is a strong Serre subcategory of $\sM$ iff the natural functor
$$ D^+(\sC) \to D_{\sC}^+(\sM)$$
is fully faithful.
\end{prop}

\noindent We recall the following theorem (not stated in its full
generality) about existence of adjoints which essentially follows
from Freyd adjoint theorem.
\begin{thm}[Freyd Adjoint Functor theorem]
Let $F:\sA \to \sB$ be any functor of abelian categories. Assume $\sA$ is cocomplete, i.e. arbitrary small colimits exists in $\sA$ and that $\sA$ has a generator. Then $F$ has a right adjoint iff it preserves colimits.
\end{thm}

\begin{rem}
In the situations of our interest the functor $F$ above will usually
be a fully faithful embedding of the category $\sC$ of
quasi-coherent sheaves on an Artin stack in a bigger category $\sM$,
where $\sM$ will either be the category of cartesian $\sO_X$-modules
on the lisse-\'etale site of $X$ or the category of all
$\sO_X$-modules on the lisse-\'etale site of $X$ (cf. \cite{lm00},
\cite{ol07}). In these cases the functor $\sC \to \sM$ is a fully
faithful embedding which preserves colimits. Moreover $\sC$ is a
cocomplete category. The existence of a generator in $\sC$ can be
proved using the fact that any quasi-coherent sheaf is a colimit of
coherent sheaves together with the observation that isomorphism
classes of coherent sheaves on a stack form a `set'. Thus, in this
situation $\sC \to \sM$ will have a right adjoint. Such embeddings
$\sC \to \sM$ will also be the main examples of strong Serre
subcategories in this paper. For more background, see \cite{lm00}
and \cite{ol07}.
\end{rem}

\begin{thm}
Let $X/\Z$ be an algebraic stack such that the diagonal
$$\Delta_{X/\Z}:  X\to X\times_{\Z}X$$ is affine. Let $\sC=QC(X)$, $\sM od^{cart}(X)$ and $\sM od(X)$ be the category of all cartesian $\sO_X$-modules and all $\sO_X$-modules, respectively,  on the lisse-\'etale site of $X$. Then $\sC$ is a strong Serre subcategory of $\sM od^{cart}(X)$ and $\sM od(X)$.
\end{thm}

\begin{proof}We first prove the case when $X$ is a affine scheme. For $\sM od^{cart}(X)$, it essentially follows from the fact that every injective quasicoherent sheaf on $X$ is flabby (cf. \cite{gr57}); for $\sM od(X)$, we observe that the inclusion of $\sM od^{cart}(X)$ in $\sM od(X)$ has a exact right adjoint. The remaining part of the proof for $\sM od^{cart}(X)$ and $\sM od(X)$ is the same, thus we will only argue for $\sM od^{cart}(X)$. Now for the general case, let $Y\to X$ be any fppf cover where $Y$ is an affine scheme. We let $Y_{\bullet}$ denote the simplicial scheme which is the $0$-coskeleton of this cover. Concretely
$$Y_i=Y\times_XY\times_X\cdots\times_X Y \ \ \ (i{\rm -times}) $$
Let $f_i:Y_i\to X$ denote the natural morphism. Note that by
assumptions on $X$, each $Y_i$ is an affine scheme and moreover the
morphisms $f_i$ are also affine. We denote by $f_{i*}$ the functor
$QC(Y_i)\to \sM od^{cart}(X)$ and $f'_{i*}$ the functor from $\sM
od^{cart}(Y_i)\to \sM od^{cart}(X)$. We claim that for any
quasi-coherent sheaf $\sG$ on $Y_i$,
\begin{equation}\label{eqn1}
R^kf_{i*}(\sG)= R^kf'_{i*}(\sG)=0 \ \ \forall \ i>0
\end{equation}
The vanishing of $R^kf_{i*}(\sG)$ for $i>0$ follows from the fact that $f_i$ is affine. The equality $$R^kf_{i*}(\sG)=R^kf'_{i*}(\sG)$$ follows from the fact that any injective sheaf on $Y_i$ is flabby. \\

\noindent For any $\sF \in QC(X)$,
there is a simplicial resolution
$$\sF \to f_{1*}f_1^*(\sF)\to f_{2*}f_2^{*}(\sF)\to \cdots,$$
Moreover, for any $\sG\in \sC$, we have the following two spectral sequences
$$Ext_{X}^p(\sG,f_{q*}f_q^*\sF)\Rightarrow Ext_{X}^{p+q}(\sG, \sF)$$
$$\tilde{Ext}_{X}^p(\sG,f_{q*}f_q^*\sF)\Rightarrow \tilde{Ext}_{X}^{p+q}(\sG, \sF)$$
where henceforth for simplicity we write
$$ Ext_X^p(-,-)=Ext^p_{QC(X)}(-,-)$$
$$\tilde{Ext}_X^p(-,-)=Ext^p_{\sM od^{cart}(X)}(-,-)$$

\noindent There is a natural morphism from the first spectral sequence to the second and to prove the theorem it is enough to show that for any $q$
$$ Ext_X^p(\sG,f_{q*}f_q^*\sF) \to \tilde{Ext}_X^p(\sG,f_{q*}f_q^*\sF)$$
is an isomorphism. But by $(1)$, we have natural isomorphisms
$$ Ext_X^p(\sG,f_{q*}f_q^*\sF) \cong Ext_{Y_q}(f_q^*\sG,f_q^*\sF)$$
$$ \tilde{Ext}_X^p(\sG,f_{q*}f_q^*\sF) \cong \tilde{Ext}_{Y_q}(f_q^*\sG,f_q^*\sF)$$
But $Y_q$ is an affine scheme. Hence we already know that $QC(Y_q)$ is a strong Serre subcategory of $\sM od^{cart}(Y)$, in particular that
$$Ext_{Y_q}(f_q^*\sG,f_q^*\sF)=\tilde{Ext}_{Y_q}(f_q^*\sG,f_q^*\sF).$$
This proves the result.
\end{proof}

\begin{proof}[Proof of Theorem (\ref{thm:sssub})] We denote by $G:\sM \to \sC$ the right adjoint of $F$. By Proposition (\ref{gsiff}) we already know that
$$RF^+: D^+(\sC) \to D_{\sC}^+(\sM)$$ is fully faithful. The assumption that $\sC$ has
enough injective objects implies $RF^+$ is essential surjective and hence an equivalence. The inverse of $RF^+$ is given by the restriction of $RG^+$ to $D_{\sC}^+(\sM)$. We now need to prove that $$RF:D(\sC) \to D_{\sC}(\sM)$$ is an equivalence for which we will apply (\ref{thm:ab4*-n}).\\

\noindent \underline{Step $1$}: \ For any object $A\in D(\sC)$,
$F(A)=0$ implies $A=0$. This implies that $F: D(\sC)\to
D_{\sC}(\sM)$ is full and essentially surjective we need to show
that the adjugant $F\circ RG \to Id$ is an equivalence. Let $D\in
D_{\sC}(\sM)$. We will show that the natural map $F(RG(D))\to D$ is
an isomorphism. We fix the following notation.
\begin{enumerate}
 \item[-] $I^{\bullet}$ a $K$-injective complex representing $D$.
 \item[-] $C=RG(D)$. Thus $C^{\bullet}=G(I^{\bullet})$ represents $C$.
 \item[-] $I_j^{\bullet}$ an injective bounded below complex representing $D_{\geq i}$.
 \item[-] For any complex $A^{\bullet}$, $^sA^{\bullet}_{\geq i}$ denotes the $i$-th stupid trunction.
 \item[-] For a complex $A^{\bullet}$, the same letter without the dot $A$, will denote the class in the derived category.
 \item[-] For a fixed $i$, we have the inverse system of complexes $\{^sG(I_j^{\bullet})_{\geq i}\}_j$. We define $M_i^{\bullet}$ to be an object such that the following triangle is a distinguished triangle.
$$ M_i^{\bullet} \to \Pi_j^{\sC} {^sG}(I_j^{\bullet})_{\geq i} \stackrel{1-{\rm shift}}{\longrightarrow} \Pi_j^{\sC}{^sG}(I_j^{\bullet})_{\geq i}$$
\end{enumerate}

\noindent \underline{Step $2$}: \ We compare $M_i$ and $C$. We claim
that there is a map $M_i \to C$ such that $H^k(M_i) \to H^k(C)$ is
an isomorphism for all $k > i$. This follows from the fact that
stupid truncation commutes with product and we have the following
diagram of distinguished triangles.
$$\xymatrix{
M_i \ar[r]\ar[d] & ^s(\Pi_j^{\sC} G(I_j))_{\geq i}\ar[r]^{1-{\rm shift}}\ar[d] & ^s(\Pi_j^{\sC} G(I_j))_{\geq i} \ar[d] \\
C \ar[r]         & \Pi_j^{\sC} G(I_j)\ar[r]^{1-{\rm shift}} &
\Pi_j^{\sC} G(I_j) }$$ The bottom triangle is a distinguished
triangle since it is obtained by applying $RG$ to the following
distinguished triangle 
$$ D \to \Pi_j I_j \stackrel{1-{\rm shift}}{\longrightarrow} \Pi_j I_j $$
which is from Theorem (\ref{thm:ab4*-n}) and
the assumption that $\sM$ is AB4*-$n$.

\noindent \underline{Step $3$}: \ We now compare $M_i^{\bullet}$ and $G(I_i^{\bullet})$. We claim that there is a map $M_i\to G(I_i)$ which induces an isomorphism
$$H^k(M_i)\to H^k(G(I_i)) \ \ \forall \ k > i+n+1 $$
Clearly for $j\leq i$, ${^sG}(I_j^{\bullet})_{\geq i} \to G(I_i)$ induces an isomorphism on the $k$-th cohomology if $k \geq i+1$. After taking prodcut over $j$, we get a map $$ \Pi_j^{\sC}({^sG}(I_j^{\bullet})_{\geq i}) \to \Pi_j^{\sC}I_i^{\bullet}$$
Since $\sC$ is AB4*-$n$, the above map induces an isomorphism
$$ H^k(\Pi_j^{\sC}{^sG}(I_j^{\bullet})_{\geq i}) \to H^k(\Pi_j^{\sC}I_i^{\bullet})  \ \ \ \ k\ge i+n+1 $$
The claim now follows from the following diagram of distinguished
triangles.
$$\xymatrix{
M_i \ar[r]\ar[d] & \Pi_{j\leq i}^{\sC} {^sG}(I_j)_{\geq i}\ar[r]^{1-{\rm shift}}\ar[d] & \Pi_{j\leq i}^{\sC} {^sG}(I_j)_{\geq i} \ar[d]\\
G(I_i) \ar[r]    & \Pi_j^{\sC}G(I_i)\ar[r]^{1-{\rm shift}} & \Pi_j^{\sC}G(I_i)\\
}$$

\noindent \underline{Step $4$}: Now consider the following commutative diagram
$$\xymatrix{
F(M_i) \ar[r]\ar[rd] & F(C)\ar[d] \ar[r] & D\ar[ld] \ar[d]\\
                     & F\circ G(I_i)\ar[r]& I_i
}$$
From Steps $2$ and $3$, and since $F$ is exact, it follows that
$$H^k(F(C)) \to H^k(F\circ G(I_i))$$ is an isomorphism for $k>i+n+1$. Moreover, since $I_i$ is a bounded below complex, by Propsition (\ref{gsiff}), $F\circ G(I_i) \to I_i$ is a quasi-isomorphism. Also, by definition of $I_i$, $H^k(D)\to H^k(I_i)$ is an isomorphism for $k\geq i$. Hence we conclude that
$$ H^k(F(C)) \to H^k(D)$$
is an isomorphism for all $k > i+n+1$. Since $i$ is arbitrary, this proves the theorem.
\end{proof}

\noindent Finally we give an example of a subcategory which is not a
strong Serre subcategory. The following example also shows that the
hypothesis that $\Delta_{X/\Z}$ is affine in Theorem
(\ref{thm:sssub}) cannot be dropped.

\begin{ex}
Let $k$ be a field and let $A/k$ be any abelian variety of positive dimension. In particular $H^1(A,\sO_X)\neq 0$. We let $X=[\Spec(k)/A]$ and $\sC$ denote the category of quasi-coherent sheaves on $X$. Let $\sM$ denote the category of all $\sO_X$-modules on the lisse-\'etale site of $X$. $\sC$ is a Serre subcategory of $\sM$. We claim that it is not a strong Serre subcategory. \\

\noindent To see this, let $\pi:\Spec(k)\to X$ denote the natural projection. The functor $\pi_*$ from lisse-\'etale sheaves on $\Spec(k)$ to those on $X$ is smooth and hence has an exact left adjoint. Therefore $\pi_*$ maps injectives to injectives and thus we have a Leray spectral sequence
$$ E_2^{p,q} = H^p(X,R^q\pi_*\sO_X) \implies H^{p+q}(\Spec(k),\sO_X)$$
which, together with the vanishing of $H^1(\Spec(k),\sO_X)$ gives us an injection
$$ 0 \to H^0(X,R^1\pi_*\sO_X) \to H^2(X,\sO_X) $$
But one can show that $H^0(X,R^1\pi_*\sO_X)=H^1(A,\sO_X)\neq 0$.
Thus $Ext^2_{\sM}(\sO_X,\sO_X)\neq 0$. However, $\sC$ itself is a
semisimple category and hence $Ext^2_{\sC}(\sO_X,\sO_X)=0$.
\end{ex}


\vspace{1cm}

\noindent Amit Hogadi\\ 
Tata Institute of Fundamental Research, Homi Bhabha Road, Colaba, Mumbai 400005. India. Email: amit@math.tifr.res.in\\
{\it Current Address}
The Mathematical Sciences Research Institute, 17 Gauss Way, Berkeley, CA 94720-5070\\

\noindent Chenyang Xu\\
Massachusetts Institute of Technology, Department of Mathematics, 77 Massachusetts Avenue, Cambridge, MA 02139-4307.  Email: cyxu@mit.edu\\
{\it Current Address}
The Mathematical Sciences Research Institute, 17 Gauss Way, Berkeley, CA 94720-5070
\end{document}